\documentclass{amsart}
\usepackage{color}

\usepackage{tipa}
\usepackage{amssymb}
\usepackage{stmaryrd}
\usepackage{graphicx}
\usepackage{amsmath,amsfonts}
\usepackage{mathtools}
\usepackage{booktabs}
\usepackage{tikz}
\usepackage{setspace}
\usepackage{float}
\usetikzlibrary{arrows,shapes,chains}
\usepackage{array}
\newcolumntype{C}[1]{>{\centering}p{#1}}

\newtheorem{theorem}{Theorem}[section]

\theoremstyle{definition}

\newtheorem{example}[theorem]{Example}

\numberwithin{equation}{section}

\begin{document}

\title[An Infinite Family of Primitive Heron Triangles]{An Infinite Family of Primitive Heron Triangles with Two Sides as Perfect Squares}

\author{Yangcheng Li}
\address{School of Mathematical Sciences, South China Normal University, Guangzhou 510631, People's Republic of China}
\email{liyc@m.scnu.edu.cn}

\thanks{This research was supported by the National Natural Science Foundation of China (Grant No. 12171163).}

\date{}


\begin{abstract}
A primitive Heron triangle is a triangle with integral sides and integral area where the greatest common divisor of the lengths of its sides is $1$. By utilizing the theory of elliptic curves, we prove that there exist infinitely many primitive Heron triangles with two sides being perfect squares. In this process, we nest one elliptic curve into another and find a surprising rational point. All the Heron triangles corresponding to this rational point are primitive. This result would imply the possible existence of infinitely many primitive Heron triangles with all three sides being perfect squares.
\end{abstract}

\maketitle

\section {Introduction}
A Heron (rational) triangle has integral (rational) side lengths and integral (rational) area. A primitive triangle is an integral triangle such that the greatest common divisor of the lengths of its sides is 1. A Pythagorean triangle means a Heron triangle whose three side lengths satisfy the Diophantine equation
\[x^2+y^2=z^2.\]

Number theory and geometry have many intrinsic connections. Many problems in number theory originate from geometry, and the methods in number theory can provide solutions to geometric problems. An ancient example is the Pythagorean theorem, which illustrates that Diophantine equations can be derived from geometric problems and this inspired the proposal of Fermat's Last Theorem. Another famous example is the congruent number problem. A positive integer is called a congruent number if it is the area of a rational right triangle. The congruent number problem is closely related to elliptic curves.

Searching for geometric objects with special forms is sometimes a difficult matter, even if they seem elementary. For example, a longstanding problem is whether there exists a Heron triangle with three rational medians. Such triangles are called perfect triangles. It is already known that there exist infinitely many Heron triangles with two rational medians (refer to \cite{Buchholz-Rathbun,Hone}). However, searching for perfect triangles remains a great challenge. Another famous open problem is whether a perfect cuboid exists. A perfect cuboid refers to a rectangular box whose edges, face diagonals, and body diagonals are all positive integers. It is already known that there exist infinitely many perfect parallelepipeds (refer to \cite{Guy,Sawyer-Reiter,Sokolowsky-Vanhooft-Volkert-Reiter}). However, it is still unknown whether a perfect cuboid exists.

Another interesting problem is whether there exists a Heron triangle with all three sides being perfect squares. This problem was first clearly proposed by Sastry \cite{Sastry} in 2001. For the case of right triangles, according to Fermat's Last Theorem, there does not exist such a right triangle. Additionally, in 2013, St\u{a}nic\u{a}, Sarkar, Gupta, Maitra and Kar \cite{Stanica-Sarkar-Gupta} prove that there does not exist such an isosceles triangle. However, they found the first primitive Heron triangle with all three sides being perfect squares: $(1853^2,4380^2,4427^2)$. They found this example by running experiments (using a ``bounded'' approach) in the GNU/Linux environment with C with GMP. Subsequently, in 2018, after searching for several months using a GMP/C$++$ parallel search program,  Rathbun \cite{Rathbun} found the second such example: $(11789^2,68104^2,68595^2).$ Rathbun conjectured that there are more such Heron triangles, but they seem difficult to locate.

In 2001, Sastry \cite{Sastry} found two primitive Heron triangles with two sides being perfect squares: $(17^2,28^2,975)$, $(29^2,37^2,1122)$. In 2013, Lagneau \cite{Lagneau} used Mathematica to search and listed $13$ primitive Heron triangles with two sides being perfect squares. In 1961, Sierpi\'nski \cite{Sierpinski1} proved that there exist infinitely many Pythagorean triangles whose two legs (not the hypotenuse) are triangular numbers $t_n$, where $t_n=n(n+1)/2$. In 1964, Sierpi\'nski \cite{Sierpinski} introduced a Pythagorean triangle (given by Zarankiewicz)
\begin{equation}
	(t_{132},t_{143},t_{164})=(8778,10296,13530),               \label{1.1}
\end{equation}
whose three sides are triangular numbers. In 2019, Peng and Zhang \cite{Peng-Zhang} proved that there exist infinitely many isosceles Heron triangles with all three sides being $k$-gonal numbers, where $k=3$ or $k\geq5$.

In this paper, we prove that there exist infinitely many primitive Heron triangles with two sides being perfect squares. In this process, we nest one elliptic curve into another and find a surprising rational point. All the Heron triangles corresponding to this rational point are primitive. This result would imply the possible existence of infinitely many primitive Heron triangles with all three sides being perfect squares.

\begin{theorem}
	Let \(A\) be a fixed positive rational number. If the elliptic curve 
	\[\mathcal{E}_{A}:~Y^2=~X^3-432(A^8-A^4+1)X-1728(A^4-2)(2A^4-1)(A^4+1)\]
	has a positive rank, then there are infinitely many scalene Heron triangles with two sides being perfect squares.
\end{theorem}

The Heron triangles in Theorem 1.1 may not be primitive.

\begin{theorem} 
Let \(A\) be a fixed positive integer. If the elliptic curve 
\[\mathcal{E}_{A}':~Y^2=~X^3 - 432(A^8 + 14A^4 + 1)X - 3456(A^{12} - 33A^8 - 33A^4 + 1)\]
has a positive rank, then there exist infinitely many primitive scalene Heron triangles with two sides being perfect squares.
\end{theorem}

\begin{theorem} 
	There exist infinitely many primitive isosceles triangles with two sides being perfect squares.
\end{theorem}

\section {Proofs of the theorems}
In this section, we will present the proofs of the theorems, and after each proof, we will provide a concrete example.

\begin{proof}[\textbf{Proof of Theorem 1.1.}]
We first consider the case of rational triangles with two sides being perfect squares. By appropriate scaling, we can obtain Heron triangles. Although the Heron triangles obtained in this way may not be primitive.

We place a rational triangle in the rectangular coordinate system in the manner shown in Figure 1. Let the three vertices of the rational triangle be \(V_{1}\), \(V_{2}\), and \(V_{3}\), respectively. And let the coordinates of its three vertices be \(V_{1}=(0,0)\), \(V_{2}=(r^{2},0)\), and \(V_{3}=(s,t)\), where $r,s$, and $t$ are all rational numbers. 
\begin{figure}[h]
	\begin{tikzpicture}
		\draw[line width=1pt, ->] (-0.7,0) -- (4,0) node[right] {$x$};
		\draw[line width=1pt, ->] (0,-0.7) -- (0,3) node[above] {$y$};
		\draw (0, 0) node[below left] {$O$};
		\draw (0, 0) node[above left] {};
		\draw (0, 0) node[below right] {$V_1$};
		\draw[line width=1pt] (0, 0) -- (2.3, 1.5) node[above] {$V_3$};
		\draw[line width=1pt] (2.3, 1.5) -- (3.5, 0) node[below] {$V_2$};
	\end{tikzpicture}
	\caption{A Heron triangle with two sides being perfect squares.}
\end{figure}
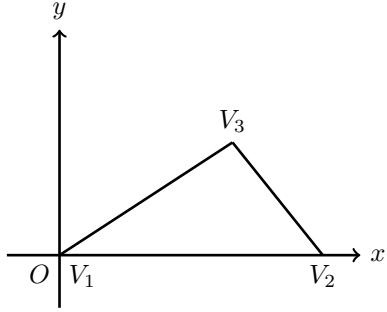

Suppose the length of side $|V_1V_3|$ is $p^2$ and the length of side $|V_2V_3|$ is $q$, where $p$ is a rational number and $q$ is a positive rational number. Then we can obtain a Diophantine equation
\begin{equation}
	\begin{cases}
		\begin{aligned}
			&s^2+t^2=p^4,\\
			&\left(s-r^2\right)^2+t^2=q^2.                                               \label{2.1}
		\end{aligned}
	\end{cases}
\end{equation}
We start by constructing infinitely many rational solutions of this system. For the first equation of (\ref{2.1}), according to the Pythagorean theorem, we obtain
\begin{equation}
s=k(m^2-1),\quad t=2km,\quad p^2=k(m^2+1),                                       \label{2.2}
\end{equation}
where $m$ is a rational number. Let $p=Ar+B$, then
\[k=\frac{p^2}{m^2+1}, \quad r=\frac{p-B}{A},\]
where \(A\) is a fixed positive rational number and $B$ is a fixed non-zero rational number. Equation (\ref{2.1}) now becomes
\begin{equation}
(m^2+1)A^4q^2=\left(A^2p^2 - (B - p)^2\right)^2m^2 + \left(A^2p^2 + (B - p)^2\right)^2,         \label{2.3}
\end{equation}
then
\[q=\frac{\sqrt{\Delta(q)}}{(m^2+1)A^2},\]
where
\[\begin{split}
	\Delta(q)=&~(m^2+1)\left(A^2p^2 - (B - p)^2\right)^2m^2 + \left(A^2p^2 + (B - p)^2\right)^2\\
	=&~(m^2 + 1)((A^2 - 1)^2m^2 + (A^2 + 1)^2)p^4\\
	&+ 4(m^2 + 1)B((A^2 - 1)m^2 -(A^2 + 1))p^3\\
	&- 2(m^2 + 1)B^2((A^2 - 3)m^2 - (A^2 + 3))p^2\\
	& - 4(m^2 + 1)^2B^3p + (m^2 + 1)^2B^4.
\end{split}\]
To get the rational solutions of equation (\ref{2.3}) with respect to $q$, it needs $\Delta(q)$ be a perfect square, say $w^2$.

Let us consider the rational points on the quartic curve
\begin{equation}
	\mathcal{C}_0: D^2=(m^2+1)\left((A^2-1)^2m^2+(A^2+1)^2\right).               \label{2.2a}
\end{equation}
By the map $\varphi_0$
\begin{align*}
	X=&~\frac{6((A^4+1)m^2+3(A^2+1)(A^2-D+1))}{m^2},\\
	Y=&-\frac{108(A^2+1)((A^4+1)m^2+(A^2+1)(A^2-D+1))}{m^3},     
\end{align*}
with the inverse map $\varphi_0^{-1}$
\begin{equation}
	\begin{split}
	m=&~\frac{6(A^2+1)Y}{(12A^4+X-24)(24A^4-X-12)},\\
	D=&~\frac{(A^2+1)(12A^4+36A^2+X+12)(12A^4-36A^2+X+12)}{(12A^4+X-24)(24A^4-X-12)}.       \label{2.3a}
\end{split}
\end{equation}
We can transform $\mathcal{C}_0$ into the elliptic curve
\begin{align*}
	\mathcal{E}_{A}:~Y^2=~X^3-432(A^8-A^4+1)X-1728(A^4-2)(2A^4-1)(A^4+1).
\end{align*}
The discriminant of $\mathcal{E}_{A}$ is
\[\Delta(A)=2176782336A^8(A^4-1)^2.\]
If $A\neq\pm1$, then $\Delta(A)\neq0$, so $\mathcal{E}_{A}$ is non-singular.

If there are infinitely many rational points on the elliptic curve $\mathcal{E}_{A}$, we can then transform the curve \(\mathcal{C}: w^2=\Delta(q)\) into an elliptic curve. Specifically, by the map $\varphi$
\begin{align*}
	X=&-\frac{3(m^2 + 1)(((A^2 - 3)m^2 - (A^2 + 3))p^2 - 3B(m^2 + 1)(B - 2p) + 3w)}{2A^2p^2},\\
	Y=&~\frac{27(m^2 + 1)^2(B - p)(((A^2 - 1)m^2 - (A^2 + 1))p^2 - B(m^2 + 1)(B - 2p) + w)}{2p^3A^3}.
\end{align*}
with the inverse map $\varphi^{-1}$
\begin{equation}
	\begin{split}
		p=&-\frac{12B(m^2 + 1)A^2(AY + 3(m^2 + 1)(X-3m^4 + 3))}{p_1(X)},\\
		w=&-\frac{16B^2(m^2 + 1)A^6w_1(X, Y)w_2(X, Y)}{(p_1(X))^2},   
	\end{split}
\end{equation}
where 
\begin{equation*}
	\begin{split}
p_1(X)=&~(2A^2X + 3(m^2 + 1)((A^2 - 3)m^2 - (A^2+3)- 3D))\\
&\times(2A^2X + 3(m^2 + 1)((A^2 - 3)m^2 - (A^2+ 3) + 3D )),\\
w_1(X)=&~X^2 - 6(m^4 - 1)X + 9(m^4 + 7m^2 + 1)(m^2 + 1)^2,\\
w_2(X, Y)=&~6A(m^2 + 1)Y + A^2X^2 + 3(m^2 + 1)((m^2 - 1)A^2 + 3m^2 + 3)X\\
 &- 9(m^2 + 1)^2((2m^2 + 1)(m^2 + 2)A^2 + 3m^4 - 3).
\end{split}
\end{equation*}

We can transform $\mathcal{C}$ into the elliptic curve
\begin{align*}
	\mathcal{E}_{m}:~Y^2=~X^3 &- 27(m^2 + 1)^2(m^4 + m^2 + 1)X\\
	& + 27(m^2 -1)(m^2 + 1)^3(2m^2 + 1)(m^2 + 2).
\end{align*}
The discriminant of $\mathcal{E}_{(A,m)}$ is
\[\Delta(m)=531441m^4(m^2 + 1)^8.\]
Hence, $\Delta(m)\neq0$, so $\mathcal{E}_{m}$ is non-singular. It is easy to verify that the elliptic curve $\mathcal{E}_{m}$ contains the following rational points:
\begin{align*}
	P_0=&\left(3m^4-3,0\right),\\
	P_1=&\left(3m^4 + 9m^2 + 6,0\right),\\
	P_2=&\left(-6m^4 - 9m^2 - 3,0\right),\\	
	Q=&\left(\frac{3(4A^2(m^4 - 1) + 3D^2)}{4A^2}, -\frac{27D(m^2 + 1)^2(A^4 - 1)}{8A^3}\right).
\end{align*}
From the rational point $Q$, we get
\begin{equation}
	p=\frac{2BD}{(m^2+1)(A^4-1)+2D}.                         \label{2.4}
\end{equation}
By some calculations, we obtain
\begin{equation}
\begin{split}
	&r=-\frac{B(m^2+1)(A^4-1)}{((m^2+1)(A^4-1)+2D)A},\\
	&q=\frac{B^2(m^2+1)((A^2-1)^4m^2+(A^2+1)^4)}{((m^2+1)(A^4-1)+2D)^2A^2}.              \label{2.5}
\end{split}
\end{equation}
Obviously, $q$ is a positive rational number. Thus, the sides of the rational triangles are $(p^2,r^2,q)$, where $p,q,r$ are given by (\ref{2.4}) and (\ref{2.5}). If the elliptic curve $\mathcal{E}_{A}$ has a positive rank, then there are infinitely many rational points on the elliptic curve $\mathcal{E}_{A}$. Each rational point on the elliptic curve $\mathcal{E}_{A}$ corresponds to a rational triangle with two sides being perfect squares. By appropriate scaling, we can obtain Heron triangles. Therefore, there exist infinitely many Heron triangles with two sides being perfect squares.
\end{proof}

\begin{example}
When $A=4$, we get an elliptic curve
\[\mathcal{E}_{4}:~Y^2=X^3-28201392X-57640996224.\]
By calculations using the package Magma, the rank of $\mathcal{E}_{4}$ is $1$ with a rational point 
\[P=(36120,-6785856).\]
From (\ref{2.3a}), we have
\[m=\frac{33}{56},\quad D=-\frac{69745}{3136}.\]
then
\[p=-\frac{2146B}{14429},\quad r=-\frac{16575B}{57716},\quad q=\frac{317052481B^2}{3331136656}.\]
By appropriate scaling, we obtain a primitive Heron triangle
\[(8584^2,16575^2,317052481)\]
where two of its sides are perfect squares. It is worth mentioning that for its three sides, any two of them are coprime to each other.
\end{example}

\begin{proof}[\textbf{Proof of Theorem 1.2.}]
In Theorem 1.1, we have obtained infinitely many Heron triangles with two sides being perfect squares. Therefore, we only need to find the primitive Heron triangles among those obtained in Theorem 1.1. In fact, the primitive Heron triangles can be obtained from (\ref{2.4}) and (\ref{2.5}).

In order to obtain the primitive Heron triangle, in (\ref{2.2a}), we assume that $A$ is a positive integer and take
\[m =\frac{t^2-s^2}{2ts},\]
where \(\gcd(s, t)=1\). Then we get
\[D^2=\frac{(s^2+t^2)^2((A^2-1)^2t^4+2s^2(A^4+6A^2+1)t^2+s^4(A^2-1)^2)}{16t^4s^4}.\]
Let us consider the rational points on the quartic curve
\begin{equation}
	\mathcal{C}_1: w^2=(A^2-1)^2t^4+2s^2(A^4+6A^2+1)t^2+s^4(A^2-1)^2.              
\end{equation}
By the map $\varphi_1$
\begin{align*}
	X=&~\frac{6((A^4 + 6A^2 + 1)t^2 + 3(A^2 - 1)^2 - 3(A^2 - 1)w)}{t^2},\\
	Y=&-\frac{108(A^2 - 1)((A^4 + 6A^2 + 1)t^2 + (A^2 - 1)^2 - (A^2 - 1)w)}{t^3},     
\end{align*}
with the inverse map $\varphi_1^{-1}$
\begin{equation}
	\begin{split}
		t=&~\frac{6s(A^2 - 1)Y}{(12A^4 - 72A^2 + X + 12)(24A^4 - X + 24)},\\
		w=&~\frac{(A^2 - 1)w_1w_2}{(12A^4 - 72A^2 + X + 12)(24A^4 - X + 24)},             \label{2.10}
	\end{split}
\end{equation}
where
\begin{equation*}
	\begin{split}
		w_1=&~12A^4 + 72A^3 + 72A^2 + 72A + X + 12,\\
		w_2=&~12A^4 - 72A^3 + 72A^2 - 72A + X + 12.
	\end{split}
\end{equation*}
We can transform $\mathcal{C}_1$ into the elliptic curve
\begin{align*}
	\mathcal{E}_{A}':~Y^2=~X^3 - 432(A^8 + 14A^4 + 1)X - 3456(A^{12} - 33A^8 - 33A^4 + 1).
\end{align*}
The discriminant of $\mathcal{E}_{A}'$ is
\[\Delta'(A)=34828517376A^4(A^4 - 1)^4.\]
If $A\neq\pm1$, then $\Delta'(A)\neq0$, so $\mathcal{E}_{A}'$ is non-singular.

From (\ref{2.4}) and (\ref{2.5}), by appropriate scaling, we obtain
\begin{equation}
	\begin{split}
	p^2=&~2^2A^2((A^2 - 1)^2(s^2 + t^2)^2 + 2^4A^2s^2t^2),\\
	r^2=&~(A^2 - 1)^2(A^2 + 1)^2(s^2 + t^2)^2,\\	
	q=&~(A^2 - 1)^4(s^2 + t^2)^2 + 2^5A^2(A^4 + 1)s^2t^2.          \label{2.11}
	\end{split}
\end{equation}

In order to obtain primitive Heron triangles, we only need to show that for any prime number $p_0$ satisfying $p_0\mid\gcd(p^{2}, r^{2})$, it always holds that $p_0^{2}\mid q$. We discuss it in the following three cases.

\begin{itemize}
	\item[$\circ$] \underline{Case~1.}
	If $p_{0}\mid A^{2} - 1$, and since $p_{0}^{2}\mid p^{2}$, we have $p_{0}^{2}\mid 2^{6}A^{4}s^{2}t^{2}$.
	
	If $p_{0} = 2$, then $p_{0}^{2}\mid q$.
	
	If $p_{0}\neq 2$, since $\gcd(A^{4}, A^{2} - 1) = 1$, it follows that $p_{0}^{2}\mid s^{2}t^{2}$. Therefore, $p_{0}^{2}\mid q$. 
	
	\item[$\circ$] \underline{Case~2.}
	Suppose $p_{0}\mid A^{2} + 1$.
	
	If $p_{0} = 2$. When $A\equiv0\pmod{2}$, $A^{2} + 1\equiv1\pmod{2}$, so $p_{0}\neq 2$. When $A \equiv1\pmod{2}$, $A^{2} + 1\equiv0\pmod{2}$, and $A^{2} - 1\equiv0\pmod{2}$, thus $p_{0}^2\mid q$.
	
	If $p_{0} > 2$, since $\gcd(A^{2}, A^{2} + 1) = 1$, we have
	\[p_{0}^{2}\mid (A^{2} - 1)^{2}(s^{2} + t^{2})^{2} + 2^{4}A^{2}s^{2}t^{2}.\]
	Because
	\[q=(A^2-1)^2\frac{p^2}{2^2A^2}+2^4A^2s^2t^2(A^2+1)^2,\]
	so $p_{0}^{2}\mid q$. 
	\item[$\circ$] \underline{Case~3.}
	If $p_{0}\mid s^{2} + t^{2}$, and since $p_{0}^{2}\mid p^{2}$, we have $p_{0}^{2}\mid 2^{6}A^{4}s^{2}t^{2}$.
	
	If $p_{0} = 2$, then $p_{0}^{2}\mid q$.
	
	If $p_{0}\neq 2$, since $\gcd(s^{2} + t^{2}, st) = 1$, it follows that $p_{0}^{2}\mid A^{4}$, so $p_{0}\mid A^{2}$. Furthermore, $p_{0}^{2}\mid A^{2}$, and thus $p_{0}^{2}\mid q$. 
\end{itemize} 

Therefore, after eliminating a certain square factor, we can obtain a primitive Heron triangle. If the elliptic curve $\mathcal{E}_{A}'$ has a positive rank, then there are infinitely many rational points on the elliptic curve $\mathcal{E}_{A}'$. Each rational point on the elliptic curve $\mathcal{E}_{A}'$ corresponds to a primitive scalene Heron triangles with two sides being perfect squares. Therefore, there exist infinitely many primitive scalene Heron triangles with two sides being perfect squares. This completes the proof.
\end{proof}

\begin{example}
	When $A=4$, we get an elliptic curve
	\[\mathcal{E}_{4}':~Y^2=X^3- 29860272X-50478615936.\]
	By calculations using the package Magma, the rank of $\mathcal{E}_{4}'$ is $1$ with a rational point 
	\[P=(6780, -242352).\]
	From (\ref{2.10}), we have
	\[t=45,\quad s=11.\]
	From (\ref{2.11}), by appropriate scaling, we obtain a primitive Heron triangle
	\[(520^2,1073^2,1020321)\]
	where two of its sides are perfect squares. It is worth mentioning that for its three sides, any two of them are coprime to each other.
\end{example}

\begin{proof}[\textbf{Proof of Theorem 1.3.}]
Suppose the three vertices of the isosceles Heron triangle are \(V_{1}=(0,0)\), \(V_{2}=(\frac{r}{2},t)\), \(V_{3}=(r,t)\). Let the side lengths \(|V_{1}V_{3}| = |V_{2}V_{3}| = p^{2}\), where \(p\) is an integer and \(r\) is a positive integer. Then we obtain the Diophantine equation
\[\left(\frac{r}{2}\right)^2+t^2=p^4.\]
According to the Pythagorean theorem, we obtain
\[p^2=m^2+n^2,\quad r=2(m^2-n^2),\quad t=2mn.\]
Take $m=a^2-1$ and $n=2a$, and thus
\begin{equation}
	p=a^2+1,\quad r=2(a^4-6a^2+1),\quad t=4(a^2-1)a,            \label{2.8}
\end{equation}
where $a$ is a positive integer and $a>1+\sqrt{2}$. If $a\equiv0\pmod{2}$, we have
\[\begin{split}
	\gcd{(p,r)}&=\gcd{(a^2+1,2(a^4-6a^2+1))}\\
	&=\gcd{(a^2+1,2(a^2+1)^2-16(a^2+1)+2^4)}\\
	&=\gcd{(a^2+1,2^4)}=1.
\end{split}\]
Therefore, $\gcd{(p^2,r)}=1.$ Hence, there exist infinitely many primitive isosceles Heron triangles whose sides are $(p^2,p^2,r)$, where $p,r$ are given by (\ref{2.8}).
\end{proof}

\begin{example}
When $a=2$, we get the primitive isosceles Heron triangles with sides $(17^2,17^2,322)$.
\end{example}

\end{document}